\documentclass[12pt]{amsart}
\usepackage{MnSymbol}
\usepackage{mathrsfs}
\usepackage[utf8]{inputenc}
\usepackage{amsmath,amsthm}

\numberwithin{equation}{section}

\begin{document}

\title{An example of a fractal finitely generated solvable group}
\author{Roman Mikhailov}

\begin{abstract}
A finitely generated solvable group with unbounded iterated identity is constructed.
\end{abstract}

\maketitle

\section{Introduction}

In \cite{Erschler}, the theory of iterated identities for groups is introduced. Recall the basic notions from \cite{Erschler}. Let $w=w(x_1,\dots, x_n)$ be a word on $n$ letters, $n\geq 1$. Consider its iterations defined as follows:
\begin{align*}
& w^{(0)}(x_1,\dots, x_n):=w(x_1,\dots, x_n),\\ & w^{(i+1)}(x_1,\dots, x_n)=w(w^{(i)}(x_1,\dots, x_n),x_2,\dots, x_n).
\end{align*}
We say that a group $G$ satisfies an iterated identity $w$, if, for any $x_1,\dots, x_n,$ there exists a number $j$ (which depends, in general, on the collection $x_1,\dots, x_n$), such that
$$
w^{(j)}(x_1,\dots, x_n)=1.
$$
In \cite{Erschler} the iterated identities are called {\it E-type (or Engel type) iterated identities}, since the iterations of the word $[x_1,x_2]$ give the Engel brackets $[[\dots [x_1,x_2],x_2],\dots,]x_2]$. 
A group $G$ is called {\it bounded}, if, for any iterated identity $w(x_1,\dots, x_n)$ of $G$, there is a number $j$ (which depends only on $G$), such that $w^{(j)}(x_1,\dots, x_n)=1$ for all collections of elements $x_1,\dots, x_n$ in $G$. A group $G$ is called {\it fractal}, if it satisfies some iterated identity and it is not bounded. The simplest example of the fractal group is the quasi-cyclic group $\mathbb Z_{p^\infty}=\varinjlim_k\ \mathbb Z/p^k$. It satisfies the iterated identity $w(x)=x^p$, but not bounded, since the orders of elements of $\mathbb Z_{p^\infty}$ are not bounded. In the same way, any finitely generated $p$-group of unbounded exponent, like Golod-Shafarevich group or Grigorchuk group, gives an example of a finitely generated fractal group.

It is shown in \cite{Erschler} that the finitely generated metabelian groups are bounded. It is natural to ask, is it true that all finitely generated solvable groups are bounded as well.

The main result of this paper is a proof that the property to be bounded can not be extended to the class of finitely generated solvable groups of class three. Let $p$ be a prime, $G$ a group given by the following presentation
$$
G=\langle x,y,t\ |\ x^t=x^p, y^{t^{-1}}=y^p, [x,y]=1, [[x,y^{jt^i}],x]=1, i,j\in \mathbb Z\rangle.
$$

\noindent{\bf Theorem 1.}
{\it The group $G$ is solvable of class three. Let $w(x_1,x_2,x_3)=[x_1,[x_2,x_3]]^p$. Then $w$ is an unbounded iterated identity for $G$.}

\vspace{.25cm}
The iterations of the identity $w$ are
\begin{align*}
& [x_1,[x_2,x_3]]^p\\
& [[x_1,[x_2,x_3]]^p,[x_2,x_3]]^p\\
& [[[x_1,[x_2,x_3]]^p,[x_2,x_3]]^p,[x_2,x_3]]^p
\end{align*}
etc. The group $G$ gives a needed example of a finitely generated solvable fractal group.

Recall that there are lot of properties which hold for the class of metabelian groups but do not hold for the class of solvable groups of class three. For example, any quotient of a finitely presented metabelian group is finitely presented. However, this is not true for solvable groups of class three \cite{Abels}. There exists a finitely presented solvable group of class three with unsolvable word problem \cite{Kharlampovich}. The main result of the present paper gives one more example of a property which differs classes of metabelian and solvable of class three groups. 

The author thanks A. Erschler for posing the problem and helpful comments and S. Ivanov for useful discussions.

\section{Proof of theorem 1}

Lets prove first that the group $G$ is solvable of class three. Denote by $X$ the subgroup of $G$ generated by elements $\langle x^{t^{-i}},\ i=0,1,\dots\rangle$ and by $Y$ the subgroup generated by elements $\langle y^{t^i},\ i=0,1,\dots\rangle$. The relations $x^t=x^p, y^{t^{-1}}=y^p$ imply that the subgroups $X$ and $Y$ are abelian and isomorphic to subgroup of rationals $\mathbb Z[\frac{1}{p}]=\{\frac{m}{p^k},\ m\in \mathbb Z,\ k\geq 0\}.$ Denote by $H$ the normal closure $\langle x,y\rangle^E$. The subgroup $H$ is generated by elements
$x^{t^{-i}}, y^{t^i},\ i=0,1,\dots$

The general element of $X$ is of the form $x^{mt^{-i}}$ for some $i,m\geq 0$, the general element of $Y$ is of the form $y^{mt^i},$ for some $i,m\geq 0$. For arbitrary $i,m,j,r,k,l\geq 0$, we will prove that the element
$$
[[x^{mt^{-i}}, y^{rt^j}], x^{lt^{-k}}]
$$
is trivial in $G$. Clearly, it is enough to consider $l=m=1$. Suppose that $k\geq i$. Then
$$
[[x^{t^{-i}}, y^{rt^j}], x^{t^{-k}}]=[[x^{t^{k-i}}, y^{rt^{j+k}}],x]^{t^{-k}}=[[x^{p^{k-i}},y^{rt^{k+j}}],x]^{t^{-k}}=1.
$$
In the same way, if $i>k$,
$$
[[x^{t^{-i}}, y^{rt^j}], x^{t^{-k}}]=[[x, y^{rt^{j+i}}], x^{t^{i-k}}]^{t^{-i}}=[[x,y^{rt^{j+i}}],x^{p^{i-k}}]^{t^{-i}}=1.
$$
That is, for all $a_1,a_2\in X,$ $b\in Y$, $[[a_1,b],a_2]=[a_1^b,a_2]=1.$ Therefore, the subgroup generated by $X$ and $Y$ is a quotient of their wreath product $X\wr Y$. The wreath product of any pair of abelian groups is metabelian by construction. Therefore, the subgroup $H$ is metabelian, $G$ is cyclic-by-metabelian and hence solvable of class at most three.

Observe that the commutator subgroup $G'$ is generated by $[H,H]$ together with elements $(x^{t^{-i}})^{p-1}, (y^{t^i})^{p-1},\ i=0,1,2,\dots$
 Since the subgroup $H$ is metabelian, the $j$th iteration of $w$ (for $j\geq 3$), can be rewritten on $G$ as
$$
[[[x_1,[x_2,x_3]]^p,[x_2,x_3]]^{p^{j-1}},[x_2,x_3],\dots, [x_2,x_3]\dots]
$$
(here we use that the element $[[x_1,[x_2,x_3]]^p,[x_2,x_3]]$ lies in $[H,H]$) Now we will prove that, for any element $u\in [H,H]$, there exists $i$, such that $u^{p^i}=1$. This will prove that $w$ is an iterated identity for $G$.

The subgroup $[H,H]$ is a normal closure of the elements $[x,y^{t^i}],\ i=1,2,\dots$ Since the group $[H,H]$ is abelian, it is enough to show that $[x,y^{t^i}]^{p^i}=1.$ We will prove it by induction on $i$. This is true for $i=0$, since $[x,y]=1$ in $G$. Suppose that
$$
[x,y^{t^i}]^{p^i}=1
$$
for a given $i$. Then
$$
1=[x,y^{t^i}]^{p^it}=[x^t,y^{t^{i+1}}]^{p^i}=[x^p,y^{t^{i+1}}]^{p^i}=[x,y^{t^{i+1}}]^{p^{i+1}}.
$$
Here we used the relation $[x,y^{t^{i+1}},x]=1$ which holds in $G$ by construction. We proved that $w$ is an iterated identity for $G$.

Now lets show that $w$ is unbounded. Fix some $i\geq 1$ and set $$x_1=x^{t^{-i}},\ x_2=t^{-1},\ x_3=y^{-t^i}.$$ Then $[x_2,x_3]=(y^{t^i})^{p-1}.$ For $j\geq 1$, the $j$th iterated identity $w$ applied to $x_1,x_2,x_3$, will give an element
$$
[[[x^{t^{-i}},(y^{t^i})^{p-1}]^p,(y^{t^i})^{p-1}]^{p^{j-1}},(y^{t^i})^{p-1},\dots, (y^{t^i})^{p-1}]=
[x^{t^{-i}},_j(y^{t^i})^{p-1}]^{p^{j}}.
$$
Here we use the standard notation for an Engel bracket $[a,_1b]:=[a,b],\ [a,_{n+1}b]=[[a,_nb],b].$

Lets rewrite the generators of $H$ as
$$
z_i:=x^{t^{-i}}, y_i:=y^{t^i},\ i=0,1,\dots
$$
The relations of $H$ are the following
\begin{align*}
& [z_i,y_j^s,z_k]=1, \\
& z_i=z_{i+1}^p, \\
& y_i=y_{i+1}^p,\\
& [z_0,y_0]=1,
\end{align*}
for all $i,j,k=0,1,2,\dots$, $s\in \mathbb Z$.
We consider the order of the element
$[z_i,_jy_i^{p-1}]
$
in $H$. Take the quotient of $H$ by the normal closure $\langle z_0,y_0\rangle^H$. The obtained group $\Gamma:=H/\langle z_0,y_0\rangle^H$ is isomorphic to the wreath product of two $p$-quasicylic groups
$$
\Gamma\simeq \mathbb Z_{p^\infty}\wr \mathbb Z_{p^\infty}.
$$
Denote by $Z$ the subgroup generated by $\{z_i\}_{i\geq 1},$ by $Y$ the subgroup generated by $\{y_i\}_{i\geq 1}$. We have $Z\simeq Y\simeq \mathbb Z_{p^\infty}$ and $\Gamma=Z\wr Y.$ We will use the multiplicative notation for elements of $Y$ and additive for elements of $Z$. We will write $u^{\frac{n}{p^i}+\mathbb Z}$ for the element $y_i^n$ of $Y$ and simply $\frac{n}{p^i}+\mathbb Z$ for the element $z_i^n$ of $Z$. The elements of the wreath product $Z\wr Y$ can be written as pairs $\left(\sum_{b\in \mathbb Z_{p^\infty}}a_bu^b, u^{b_0}\right),\ b_0\in \mathbb Z_{p^\infty}$ with the product given as
$$
\left(\sum_{b\in \mathbb Z_{p^\infty}}a_bu^b, u^{b_0}\right)\left(\sum_{b\in \mathbb Z_{p^\infty}}a_b'u^b, u^{b_0'}\right)=\left(\sum_{b\in \mathbb Z_{p^\infty}}a_bu^b+\sum_{b\in \mathbb Z_{p^\infty}}a_b'u^{b+b_0},u^{b_0+b_0'} \right)
$$
In this notation, the element $[z_i,_jy_i^{p-1}]$ corresponds to the expression
$$
\left((\frac{1}{p^i}+\mathbb Z)(u^{\frac{p-1}{p^i}+\mathbb Z}-1)^j,1\right).
$$
For $j<i$, the order of this element is $p^i$ [after opening the bracket, there will be no term which can cancel $\frac{(-1)^j}{p^i}$]. This implies that, for $j<i$, the order of $[z_i,_jy_i^{p-1}]$ is $p^i$.
We conclude that
$$
[[[x^{t^{-i}},(y^{t^i})^{p-1}]^p,(y^{t^i})^{p-1}]^{p^{j-1}},(y^{t^i})^{p-1},\dots, (y^{t^i})^{p-1}]\neq 1
$$
for $j<i$. That is, we need $i$ iterations of the word $w$ to become trivial for  $x_1=x^{t^{-i}},\ x_2=t^{-1},\ x_3=y^{-t^i}$. The statement is proved.


\begin{thebibliography}{99}
\bibitem{Abels} H. Abels: An example of a finitely presented solvable group, Homological group theory (Proc. Sympos.,
Durham, 1977), London Math. Soc. Lecture Note Ser., vol. 36, Cambridge Univ. Press, Cambridge, 1979,
pp. 205–211.

\bibitem{Erschler}  A. Erschler: Iterated identities and iterational depth of groups, {\it J. Mod. Dyn.} {\bf 9} (2015), 257--284.

\bibitem{Kharlampovich} O. Kharlampovich: A finitely presented solvable group with unsolvable word problem, {\it Izv. Akad. Nauk. SSSR} {\bf 45} (1981), 852--873. 
\end{thebibliography}
\end{document}